\documentstyle [theorem] {article}

\newcommand{\Complex}{
        \mbox{C \hspace{-1.16em} \raisebox{-0.018em}{\sf l}}\;}
\title{Phase Space Reduction for Star-Products:\\
        An Explicit Construction for $\Complex \mbox{\LARGE P}^n$ \vspace{1cm}}
\author{{\bf M. Bordemann\thanks{mbor@phyq1.physik.uni-freiburg.de}~,
            \addtocounter{footnote}{2}
            M. Brischle\thanks{brischle@phyq1.physik.uni-freiburg.de}~,
            \addtocounter{footnote}{-2}
             C. Emmrich\thanks{cemm@phyq1.physik.uni-freiburg.de}~,
             \addtocounter{footnote}{1}
             S. Waldmann\thanks{waldman@phyq1.physik.uni-freiburg.de}}\\[3mm]
             Fakult\"at f\"ur Physik\\Universit\"at Freiburg\\
          Hermann-Herder-Str. 3\\79104 Freiburg i.~Br., F.~R.~G\\[3mm]
       }
\date{FR-THEP-95/6 \\[1mm] March 1995}

\input {amssym.def}
\input {amssym}

\newcommand {\EQREF} [1] {eqn. (\ref{#1})}

\newcommand {\BEQ} [1] {\[}
\newcommand {\EEQ} {\]}
\newcommand {\BEQN} [1] {\begin {equation} \label {#1}}
\newcommand {\EEQN}{\end {equation}}

\newcommand {\QED} {\hfill $\blacksquare$}
\newcommand {\comment}[1]{}

\newcommand {\CP} {\Bbb{C}\mbox{P}}
\newcommand {\CNull} [1] {\Bbb{C}^{#1} \setminus \{ 0 \} }
\newcommand {\Cinfty}{C^{\infty}(\CNull{n+1})}

\theoremheaderfont {\bf}
\theoremstyle {plain}

\theorembodyfont {\it}
\newtheorem {DEF} {Def.:}[section]

\theorembodyfont {\it}
\newtheorem {THEOREM} [DEF]  {Theorem:}
\newcommand {\THEOREMREF}[1] {theorem (\ref{#1})}

\theorembodyfont {\it}
\newtheorem {PROP}  [DEF]  {Proposition:}
\newcommand {\PROPREF} [1] {proposition (\ref{#1})}

\theorembodyfont {\it}
\newtheorem {LEMMA} [DEF] {Lemma:}

\theorembodyfont {\rm}

\newenvironment {PROOF}{\small {\bf Proof:}}{\QED \\}

\newcommand {\newstar} {\,\tilde{\ast}\,}
\newcommand {\mustar} {\,\ast^{\mu}\,}

\begin {document}

\maketitle
\begin{abstract}
We derive a closed formula for a star-product on complex projective space
and on the domain $SU(n+1)/S(U(1)\times U(n))$
using a completely elementary construction:
Starting from the standard star-product of Wick type on $\CNull{n+1}$ and
performing a quantum analogue of Marsden-Weinstein reduction, we can give
an easy algebraic description of this star-product. Moreover,
going over to a modified star-product on $\CNull{n+1}$,
obtained by an equivalence transformation, this description
can be even further simplified, allowing the explicit computation of
a closed formula for the star-product on $\CP^n$ which can easily
transferred to the domain $SU(n+1)/S(U(1)\times U(n))$.
\end{abstract}
\vfill
\newpage

% ********** functions *****************

\section {Introduction}

The concept of deformation quantization has been defined by
F. Bayen, M. Flato, C. Fronsdal, A. Lichnerowicz, D. Sternheimer in 1978
(cf. \cite{bayen}) and in the more restricted context of K\"ahler
manifolds
by F. Berezin already in 1974 (\cite{Ber 74}). The basic idea is to
formally deform the pointwise commutative multiplication in the
space of smooth complex-valued functions on a symplectic manifold $M$,
$C^{\infty}(M)$,
to a noncommutative associative multiplication, the star-product,
whose first order commutator is proportional to the Poisson bracket.
In addition, the definition of a star-product usually
requires that the deformation series be local, i.e. the existence
of an atlas in which each term in the series is locally a bidifferential
operator, that the star-product with constant functions reduce to
the pointwise multiplication, and that the complex conjugation be
an antilinear involution in the deformed algebra, the most prominent
cases being of Weyl-Moyal type where the formal parameter is considered
to be purely imaginary and of Wick type where the formal parameter
is taken to be real. The Wick-type case typically occurs on  K\"ahler
manifolds.

The nontrivial question of existence of star-products on every
symplectic manifold has been positively answered by M. DeWilde and
P.B.A. Lecomte in 1983 \cite{DL 83}
using the computation of Hochschild and
Chevalley-Eilenberg cohomology of $C^{\infty}(M)$ (e.g. \cite{CDG 80})
and independently by B. Fedosov in 1985 (\cite{Fed 85} and \cite{Fed 94})
who did not use these cohomological computations.

However, closed formulas
for the deformation series as in the case of the Weyl-Moyal product
(resp. the Wick product)
on a symplectic (resp. K\"ahler) manifold admitting a flat torsion-free
symplectic (K\"ahler) connection (e.g \cite{bayen})
are commonly regarded as practically
nonexisting for all other manifolds with the following notable exceptions:
the article of S. Gutt
(\cite{Gutt 83}) for the cotangent bundle of a connected Lie group,
the work of C. Moreno and
P. Ortega-Navarro (\cite{MO 83}, \cite{Mor 86}, \cite{Mor 87}) who gave
recursion
formulas in the case of the 2-sphere, the Poincar\'e disk and their
higher dimensional analogues, complex projective space and the symmetric
hermitean domain $SU(1,n)/S(U(1)\times U(n))$ (e.g. \cite{Hel 78},
p.518 for details), and the work of M. Cahen, S. Gutt, J. Rawnsley
(\cite{CGR III})
who gave an explicit description of the star-product of Berezin symbols
in the case of the Poincar\'e disk.

Another problem which does not seem to have been much attacked in the
framework
of star products is the question whether and how the procedure
of symplectic reduction---which has turned out to be a powerful tool
in classical mechanics for both the construction of new symplectic manifolds
and Liouville integrable systems---can be extended to a ``reduction'' of
star products. For certain $U(1)$-actions this has been formulated in a
rather general, but not too explicit fashion again by B. Fedosov
(\cite{Fed 94a}). In the framework of
geometric quantization of K\"ahler manifolds and coherent states
(see the series \cite{CGR I}, \cite{CGR II}, \cite{CGR III} for details)
the fact that the K\"ahler manifold is a $U(1)$ reduction of the total
space (minus the zero-section) of its prequantum line bundle has been
exploited in \cite{BMR 93} to reformulate the asymptotic limits
already described in \cite{CGR II}.

In this paper we give an explicit closed formula for a star-product on
complex projective space and the domain $SU(1,n)/S(U(1)\times U(n))$
which in addition is shown to be the Marsden-Weinstein reduction of a
star-product equivalent to the usual
Wick product on the flat K\"ahler manifold $\CNull{n+1}$. The construction
is completely elementary and does not use geometric quantization or
deep cohomological results in complex differential geometry.

The paper is organized as follows: we show in section 2. that
the usual Wick product of two functions on $\CNull{n+1}$ which are
pulled back from complex projective space is a power series in pulled-back
functions multiplied by ``radial'' functions (which only
depend on the Euclidean distance from the origin). Moreover,
the Wick product of two radial functions is radial but not simply
pointwise multiplication. In section 3. we explicitly construct an
equivalence transformation of the Wick product on $\CNull{n+1}$ thus
obtaining a new star-product in which
the radial functions behave like scalars in the subalgebra of
$U(1)$-invariant functions. In section
4. we show that this new star-product on $\CNull{n+1}$ can simply
be projected to complex projective space by essentially fixing the square
of the
Euclidean distance (which is---up to a constant factor---nothing but the
$U(1)$ momentum map). The deformation series of this star-product on
$\CP^n$ can explicitly be written down (see Theorem 3.1). In Section 5. we
show that all the results for $\CP^n$ can easily be transferred to
the noncompact domain $SU(1,n)/S(U(1)\times U(n))$. Finally, Section 6.
contains a proof
that our deformation series solves the recursion equations of the one
implicitly given by Moreno in (\cite{MO 83}) and thus coincides with it.

\section {Special properties of the Wick product}

Let us first recall some well-known facts about the Wick-product on
$\CNull{n+1}$:

\noindent We denote by $z^0=q^0+{\bf i}p^0, \ldots, z^n
=q^n+{\bf i}p^n$ the standard co-ordinates
in $\CNull{n+1}$ and
$\partial/\partial z^k=(1/2)(\partial/\partial q^k
                       -{\bf i}\partial/\partial p^k)$
and
$\partial/\partial \bar{z}^k=(1/2)(\partial/\partial q^k
                           +{\bf i}\partial/\partial p^k)$
the corresponding complex differential operators. We shall use the
summation convention, i.e the summation over repeated indices from zero
to $n$ is automatic.
Furthermore, in this paper the symbol ``$z$'' will always designate
a point in
$\CNull{n+1}$ and for a smooth complex-valued function $F$ on
$\CNull{n+1}$ the notation $F(z)$ will {\em not} automatically
imply that $F$ is holomorphic unlike the traditional use in function
theory.
The standard symplectic form $\omega_0$ on $\CNull{n+1}$ reads
$\frac{{\bf i}}{2}dz^k\wedge d\bar{z}^k$. We shall denote by $\Cinfty$
the space of smooth complex-valued functions on $\CNull{n+1}$. The Poisson
bracket
of two functions $F$ and $G$ in $\Cinfty$ is given by
\BEQ {Poissonklammeroben}
   \{F,G\} =
   \frac{2}{{\bf i}}
  (\frac{\partial F}{\partial z^k}\frac{\partial G}{\partial \bar{z}^k}
   -\frac{\partial F}{\partial \bar{z}^k}\frac{\partial G}{\partial z^k})
\EEQ
The {\em Wick product} is defined for two
functions $F,G\in\Cinfty$ by the following formal series with formal
parameter $\lambda$:
\BEQN {WickDef}
    F \ast G := \sum_{r=0}^\infty \frac{\lambda^r}{r!}
                \frac{\partial^r F}
                {\partial z^{i_1} \ldots \partial z^{i_r}}
                \frac{\partial^r G}
                {\partial \bar{z}^{i_1} \ldots \partial \bar{z}^{i_r}}.
\EEQN
It is well-known (cf. \cite{Ber 74}) that this definition naturally extends
to an
{\em associative} product on the space of formal power
series in $\lambda$ having coefficients in $\Cinfty$ which we shall call
$\cal A$.
This sum converges e.g. for complex-valued polynomials in which case we can
identify $\lambda$ with twice Planck's constant $\hbar$ yielding the first
order commutator
\BEQ {WickCom}
    F \ast G - G \ast F = \frac{{\bf i}\lambda}{2} \{F, G \} + \ldots
                        = {\bf i}\hbar \{F, G \} + \ldots
\EEQ

In the following we consider functions $\Cinfty$ of
a particular form:

\noindent Define $x : \CNull{n+1} \to \Bbb{R}^+$ by
\BEQ {xdef}
     x(z) := \bar{z}^k z^k
\EEQ
{}From time to time we shall consider $x$ as co-ordinate function along
$\Bbb{R}^+$ on $\CNull{n+1}=\Bbb{R}^+\times S^{2n+1}$.
A function $R\in\Cinfty$ is called {\em radial}
iff there exists a $C^\infty$-function $\varrho : \Bbb{R}^+ \to \Bbb{C}$
such that $R (z) = \varrho \circ x (z)$.
The flow of the vector field $x\frac{\partial}{\partial x}
                  =\frac{1}{2}(z^k\frac{\partial}{\partial z^k}+
                               \bar{z}^k\frac{\partial}{\partial \bar{z}^k})$
consists of the scalar multiplication of a complex vector in $\CNull{n+1}$
by a positive real number. Furthermore, the flow of the vector field
$Y(z)={\bf i}(z^k\frac{\partial}{\partial z^k}-
                               \bar{z}^k\frac{\partial}{\partial \bar{z}^k})$
generates scalar multiplication of $z$ by an element of the unit circle
$U(1)$. We call a function $f$ in $\Cinfty$ {\em
homogeneous} iff it is invariant by the natural action of the group
$\Bbb{C}\setminus \{0\}$ on $\CNull{n+1}$.

By direct calculation we find the following lemma:
\begin {LEMMA} \label {specialstar}
    Let $R_1 = \varrho_2 \circ x$, $R_2 = \varrho_2 \circ x$,
    be radial functions, $F$ be a $U(1)$-invariant function in
    $\Cinfty$, and $f,g$ two homogeneous
    functions on $\CNull{n+1}$.
    Then the following relations for the
    Wick-product hold:
    \begin {enumerate}
    \item
    \BEQN {radu1}
        R_1 \ast F = \sum_{r=0}^\infty
        \frac{\lambda^r}{r!} x^r \frac{\partial^r\varrho_1}{\partial x^r}(x)
        \frac{\partial^r F}{\partial x^r}= F\ast R_1
    \EEQN
    \item In particular we have
    \BEQN {radrad}
        R_1 \ast R_2 = R_2 \ast R_1,
    \EEQN
        and $R_1 \ast R_2$ is again radial.
    \item The Wick star-product of a radial with a homogeneous function
      reduces to pointwise multiplication:
    \BEQN {radhomo}
        R_1 \ast f = R_1f = f \ast R_1.
    \EEQN
    \item For each positive integer $r\geq 0$ denote by $M_r$ the following
      bidifferential operator defined on a pair of functions $G$ and $H$ in
      $\Cinfty$:
    \BEQN {dieMs}
        M_r(G,H)
            =x^r\frac{\partial^r G}{\partial z^{i_1}\cdots \partial z^{i_r}}
   \frac{\partial^r H}{\partial \bar{z}^{i_1}\cdots \partial \bar{z}^{i_r}}.
    \EEQN
    Then each function $M_r(f,g)$ is homogeneous and we have the
    obvious formula
    \BEQN {WickMr}
        f \ast g  = \sum_{r=0}^\infty \frac{1}{r!}
                           \frac{\lambda^r}{x^r} M_r (f, g).
    \EEQN
    \end {enumerate}
\end {LEMMA}
\begin{PROOF}
Since radial and homogeneous functions are $U(1)$-invariant, and
$\partial/\partial x$ vanishes on homogeneous functions only the
first and the last part of this lemma need to be proved: denoting
by $E$ and $\bar{E}$
the Euler operators $z^i\frac{\partial}{\partial z^i}$
and $\bar{z}^i\frac{\partial}{\partial \bar{z}^i}$ we get by induction
\[     z^{i_1}\ldots z^{i_r}
       \frac{\partial^r F}{\partial z^{i_1}\ldots\partial z^{i_r}}(z)
        =  (\prod_{k=0}^{r-1}(E-k)F)(z) \\
         = (\prod_{k=0}^{r-1}(x\frac{\partial}{\partial x}-k)F)(z)
                   = (x^r\frac{\partial^r F}{\partial x^r})(z)    \]
where we have used that $x\frac{\partial}{\partial x} = (1/2)(E+ \bar{E})$
and that $Y={\bf i}(E-\bar{E})$ vanishes on $U(1)$-invariant smooth
functions. In an analogous way the corresponding statement involving the
$\bar{z}^i$ and $\partial/\partial \bar{z}^i$ is shown. This proves
part i). For the part iv), note that a
smooth complex-valued function is homogeneous iff it is annihilated by
the two Euler operators $E$ and $\bar{E}$. Clearly,
$[E,\partial/\partial z^k]=-\partial/\partial z^k$ and
$[\bar{E},\partial/\partial \bar{z}^k]=-\partial/\partial \bar{z}^k$ whereas
$[E,z^k]=z^k$ and $[\bar{E},\bar{z}^k]=\bar{z}^k$ and the homogeneity
of $M_r(f,g)$ follows by induction.
\end{PROOF}

Consider now the complex projective space $\CP^n$ with its natural projection
$\pi:\CNull{n+1}\rightarrow \CP^n:z\mapsto [z]$ where $[z]$ denotes the
complex line defined by $z$. $\CP^n$ is a K\"ahler manifold with a
complex structure inherited by $\CNull{n+1}$ and a symplectic form
$\omega$ which can be defined by its pull-back $\pi^*\omega=
\frac{{\bf i}}{2}\partial\bar{\partial}\log x$
(see e.g. \cite{GH 78} for details)
or by phase space reduction, see section 4.
The pull-back of a smooth complex-valued function $\phi$ to $\CNull{n+1}$,
$\pi^*\phi$ is obviously homogeneous, and vice versa.

Now looking at the third and forth statement of the above lemma one
might be tempted to directly project the above Wick star-product of two
homogeneous functions to a star-product on functions on complex projective
space by using $\lambda/x$ in some sense as a new formal parameter.
However, by the
first statement of the above lemma it quickly becomes clear that things
are not that simple: the Wick product of two radial functions is {\em not}
pointwise multiplication, although the Wick product of a radial and a
homogeneous function is pointwise. This fact destroys the associativity
of such a naively defined ``star-product''. But we shall see in the next
section
that an equivalence transformation modifying the Wick star-product
a little will lead out of this difficulty.

% ********** Transformation ************

\section {Equivalence transformation of the Wick-product}

The aim of this section is to find an equivalence transformation $S$ (see
e.g. \cite{bayen}, part I, p.85 for definitions) of
the Wick-product $\ast$ into a new star-product $\newstar$ with
the property that the $\newstar$-product of two radial functions is
the ordinary pointwise product. By \EQREF{radu1} und
\EQREF{radrad} we only have
to consider the associative local deformation of the pointwise product of
two smooth complex-valued functions
$\varrho_1,\varrho_2$ on the positive real line defined by
\BEQN{Radstar}
     \varrho_1\star\varrho_2(x):=\sum_{r=0}^{\infty}
          \lambda^r \frac{x^r}{r!}
                \frac{\partial^r\varrho_1}{\partial x^r}(x)
                \frac{\partial^r\varrho_2}{\partial x^r}(x).
\EEQN
For general reasons, (cf. e.g. \cite [p. 46] {gerstenhaber}) this
deformation must be equivalent to the pointwise multiplication
because the second local
Hochschild cohomology group of the associative algebra of smooth
complex-valued functions on $\Bbb{R}^+$ vanishes being isomorphic to
$\Gamma(\bigwedge^2 T \Bbb{R}^+) = \{0\}$ (cf. e.g. \cite{CDG 80}).

Such an equivalence transformation $S$ is a formal power series
$S(x, \partial_x) =\sum_{r=0}^\infty \lambda^r S_r(x,\partial_x)$
(where the $S_r$ are differential operators acting on
the space of all smooth complex-valued functions
on the positive real line and we have used the notation
$\partial_x:=\partial/\partial x$) with the property
\BEQN {equiv}
    S (\varrho_1 \star \varrho_2) = (S \varrho_1) (S \varrho_2).
\EEQN
for arbitrary smooth functions
$\varrho_1, \varrho_2: \Bbb{R}^+ \to \Bbb{C}$
and $S_0$ is equal to the identity.
In order to perform concrete calculations it will turn out to be convenient
to
work with the following {\em symbol} $\hat{S}$ of $S$ which will be
a power series
in $\lambda$ whose coefficients are contained in the space of smooth
complex-valued functions on $\Bbb{R}^+\times\Bbb{R}$ that are polynomial
in the second variable and which is defined by
\BEQ{thesymbol}
  \hat{S} (x,\alpha) e_\alpha (x) := (S e_\alpha)(x)
\EEQ
where $e_\alpha$ denotes the exponential function
$x\mapsto e^{\alpha x}$ on the real line for $\alpha\in\Bbb{R}$.
This makes sense since
$\partial_x^r e_\alpha=\alpha^r e_\alpha~~\forall r\in\Bbb{N}$. Moreover,
since $S$ is uniquely determined by its action on the monomials $x^r$ it
can obviously be regained
from its symbol $\hat{S}$ by putting the monomials in the second
variable $\alpha$ occurring in $\hat{S}(x,\alpha)$ on the right and then
substituting $\alpha$ by $\partial_x$ (the standard ordering prescription).

We get the following theorem:

\begin {THEOREM}
    The deformation $\star$ is equivalent to the pointwise product and
    every equivalence transformation $S$ satisfying \EQREF{equiv}
    has a symbol of the form
    \BEQN {Ssymbol}
        \hat {S} (x, \alpha) =
        \exp\left(\frac{x}{\lambda}\left(D(x,\lambda)
        \log (1 + \lambda\alpha) - \lambda\alpha\right) \right)
    \EEQN
    where $D(x,\lambda) := exp(\lambda C(x,\lambda))$ is an arbitrary
    formal power series in $\lambda$ with smooth coefficient functions
    starting with $1$. Moreover $S1=1$ for the constant function $1$.
\end {THEOREM}
\begin{PROOF} First note that in terms of exponential functions $e_\alpha$,
$e_\beta$ with
$\alpha,\beta\in\Bbb{R}$ the above star-product \EQREF{Radstar} on the
positive real line has the following particularly simple form:
\BEQ {e-func}
    e_\alpha \star e_\beta = e_{\alpha+\beta+\lambda\alpha\beta}
\EEQ
where the r.h.s. is well-defined as a formal series in $\lambda$. Using
this equation the condition \EQREF{equiv} can be reformulated into
the following functional equation for the symbol of $S$:
\BEQN {Behauptung}
    \hat{S} (x, \lambda\alpha\beta+\alpha+\beta) e^{\lambda\alpha\beta x}
    =
    \hat{S} (x,\alpha) \hat {S} (x,\beta)
\EEQN
Since $\hat{S}_0=1$ we may take the formal logarithm
$\lambda T = \log \hat {S}$ on both sides thus getting the following
equation for $T$,
which again is a formal power series starting with a constant:
\[
    \lambda\alpha\beta x + \lambda
    T(x, \lambda\alpha\beta + \alpha + \beta) =
    \lambda T(x, \alpha) + \lambda T(x, \beta)
\]
By differentiating this equation with respect to
$\alpha$ and $\beta$ we see that every solution of this
functional equation satisfies the following differential
equation, where $'$ denotes the differentiation with respect to the
second argument.
\[
    x + T''(x, \lambda\alpha\beta + \alpha + \beta) (1 + \lambda
    (\lambda\alpha\beta + \alpha + \beta)) + \lambda T'(x,
    \lambda\alpha\beta + \alpha + \beta) = 0
\]
With $\gamma := \lambda\alpha\beta + \alpha + \beta$ this differential
equation has the general solution
\BEQ {diffsolution}
    T(x,\gamma) = \frac{x}{\lambda^2} e^{\lambda C(x,\lambda)}
    \log (1 + \lambda\gamma) - \frac{x \gamma}{\lambda} +
    \frac{x B(x,\lambda)}{\lambda}
\EEQ
depending on two arbitrary integration constants $B(x,\lambda)$ and
$C(x,\lambda)$.
Now we can easily check that this general solution of the differential
equation leads to a solution $\hat {S}$ of the functional equation
(\ref {Behauptung}) iff $B$ vanishes identically. In this case we have
$\hat{S}_0=1$. The formal power series
$C(x,\lambda)=\sum_{r=0}^\infty \lambda^r C_r(x)$
is completely arbitrary. Since $(S1)(x)$ is obviously equal to $\hat{S}(x,0)$
which in turn is equal to $1$ the theorem is proved.
\end{PROOF}

In the next step we regard the differential operator $\partial_x$
as a differential operator on
$\Cinfty\cong C^{\infty}(\Bbb{R}^+ \times S^{2n+1})$ by observing that
$\partial_x =\frac{1}{2x}(z^k\frac{\partial}{\partial z^k}+
                            \bar{z}^k\frac{\partial}{\partial \bar{z}^k})$.
Then $S$ extends to a
formal power series in $\lambda$ with differential operators acting
in $\Cinfty$ which we shall also denote by $S$. We can now use
$S$ to transform the Wick-product $\ast$ into a new
star-product $\newstar$ in $\cal A$ defined by
\BEQN {newstar}
    F \newstar G := S \left( \left(
    S^{-1} F \right) \ast
    \left( S^{-1} G \right) \right)
\EEQN
for two functions $F, G\in {\cal A}$.
If $f = \phi \circ \pi$ is a homogeneous function we
have $\partial_x f = 0$. This leads to
\BEQ {Shomo}
    S f = S^{-1} f = f
\EEQ
Combining this result with the properties of the Wick-product in
lemma (\ref{radu1}), (\ref{radhomo}) and (\ref{WickMr}) we
get four relations for the new star-product analogous to those mentioned
in the lemma in the previous section:

\begin {THEOREM} \label {fundamental}
    Let $R_1 = \varrho_1 \circ x$, $R_2 = \varrho_2 \circ x$ be two
    radial functions, $f = \phi \circ \pi$, $g = \psi \circ \pi$
    two homogeneous functions and $F$ a $U(1)$-invariant function in
    $\Cinfty$. Then
    \BEQ {RnewstarR}
        R_1 \newstar R_2 = R_1R_2 = R_2 \newstar R_1
    \EEQ
    \BEQN {Rnewstaru1}
        R_1 \newstar F =R_1F= F \newstar R_1
    \EEQN
    \BEQ {RnewstarHomo}
        R_1 \newstar f = R_1f = f \newstar R_1
    \EEQ
    \BEQ {HomonewstarHomo}
        f \newstar g = S(f \ast g) =
         \sum_{r=0}^\infty \frac{1}{r!}
        \left(S(x,\partial_x)\left(\frac{\lambda}{x}\right)^r\right)
               M_r (f,g) .
    \EEQ
\end {THEOREM}

Because $M_r(f,g)$ is again homogeneous, the
differential operators in $S$ act on the radial functions
$(\frac{\lambda}{x})^r$ only. We can give an explicit form of these terms.
First we use the equation (\ref{Ssymbol}) for $S$ to find the
following proposition.

\begin {PROP} \label {Sxpower}
    For any formal power series $D(x, \lambda)$ starting with $1$ the
    following equations hold for the differential operator
    $S(x,\partial_x)$ defined by its symbol as in \EQREF{Ssymbol}:
    \BEQN {SxDx}
        S(x, \partial_x) x = D(x,\lambda)x
    \EEQN
    For $r \in \Bbb{N}$ we have
    \BEQ {Sxposr}
        S(x, \partial_x) x^r = (D(x,\lambda)x)^r \prod_{k=0}^{r-1}
        \left(1- k\frac{\lambda}{D(x,\lambda)x} \right)
    \EEQ
    \BEQN {Sxnegr}
        S(x,\! \partial_x) \frac{1}{x^r} \!=\!
           \left(\!\frac{1}{D(x,\lambda)x}\right)^{\! r} \!
        \prod_{k=0}^{r}\! \left(\!1 \!
        +\! k \frac{\lambda}{D(x,\lambda)x}\! \right)^{\!-1}\!\! =\!\!
           \left(\!\frac{1}{D(x,\lambda)x}\!\right)^{\!r}\!
    \sum_{s=0}^{\infty}\!\left(\!\frac{\lambda}{D(x,\lambda)x}\!\right)^{\!s}
                \!          A_s^{(r)}
    \EEQN
    where the $A_s^{(r)}$ are rational numbers defined by $A^{(0)}_0=1$,
     $A^{(0)}_s=0~~\forall s\geq 1$, and for $s\geq 0,~r\geq 1$
    \BEQN {Arsdef}
        A^{(r)}_s := \frac{1}{(r-1)!} \sum_{k=1}^r {r-1  \choose k-1}
                k^{s+r-1} (-1)^{r+s-k}.
    \EEQN
\end {PROP}
\begin {PROOF}
The first equation is proved by observing that $S(x,\partial_x)x=
(\partial/\partial \alpha)(Se_\alpha)(x)|_{\alpha=0}$ and using the
symbol $\hat{S}$.
For the second and the third statement we
use the fact that the Wick-product of $x$ and $x^r$ is
$x^{r+1}+\lambda rx^r$. We use \EQREF{newstar} to obtain
a recursion formula for $Sx^r$ and $Sx^{-r}$. These recursion formulas
can be solved using \EQREF{SxDx}. The second equation in \EQREF{Sxnegr}
is proved by performing a partial fraction decomposition of the product
and induction.
\end {PROOF}

We want to choose the formal series $D$ in such a way
that in the new star-product of two homogeneous functions
of $\lambda$ only the combination $(\frac{\lambda}{x})^r$ times
some homogeneous term depending on $f,g$ occurs in each order.
This is the case iff the series $D$ is of the form
\BEQN {theD}
    D(x, \lambda ) = \sum_{r=0}^\infty \left(
                     \frac{\lambda}{x}\right)^r d_r
                     \qquad \mbox{where } d_0 = 1
                     \mbox { and } d_r \in \Bbb{C}.
\EEQN

We shall henceforth assume that $D$ is of this special form with
arbitrary complex numbers $d_r$. Then we
get the following formula for the new star-product of two homogeneous
functions $f$ and $g$:
\BEQN {thenewstar}
    f \newstar g = \sum_{r=0}^\infty \frac{1}{r!}
                   \left(\frac{\lambda}{D(x,\lambda)x}\right)^r
                   \prod_{k=1}^r \left( 1 + k \frac{\lambda}{D(x,\lambda)x}
                   \right)^{-1} M_r (f, g)
\EEQN
where the empty product in the case $r=0$ is defined to be $1$.
This can be rewritten as
\BEQN {homohomo}
    f \newstar g = \sum_{r=0}^\infty \left(\frac{\lambda}{x}\right)^r
                   K_r (f,g)
\EEQN
where the $K_r (f,g)$ are some linear combination --- depending on
the choice of the series $D$ --- of the $M_s (f, g)$ with $s \le r$.
In particular $K_r (f,g)$ is again homogeneous.

An interesting property of the original Wick product is the fact that
it can be expressed as a power series in a differential operator
on smooth complex-valued functions $F$ on
$\CNull{n+1}\times\CNull{n+1}$
namely the operator $\cal P$ defined by
\BEQ{Wickdiff}
    {\cal P}F:=\frac{\partial^2F}{\partial z^i\partial \bar{w}^i}.
\EEQ
Upon writing $m$ for the evaluation of a smooth complex-valued
function $F$ on $\CNull{n+1}\times\CNull{n+1}$ on the diagonal (i.e
$m(F)(z):=F(z,w)|_{z=w}$) we get the following formula for the
old Wick product of $G$ and $H$ in $\cal A$:
\BEQ{Wickkurz}
      G\ast H=m e^{\lambda{\cal P}}G\otimes H           .
\EEQ
A similar thing can be done for the new star-product:
We define for a smooth complex-valued function $F$ on
$\CNull{n+1}\times\CNull{n+1}$ the following differential operators:
\BEQN {NOp}
    N(F) (z, w) := z^i\bar{w}^i
                     \frac{\partial^2 F}
                            {\partial z^j \partial \bar{w}^j}(z,w)
\EEQN
and
\BEQN {calMR}
    {\cal M}_r (F) (z, w) :=
    z^{i_1} \ldots z^{i_r}
    \bar{w}^{i_1} \ldots \bar{w}^{i_r}
    \frac{\partial^{2r} F}
    {\partial z^{j_1} \ldots \partial z^{j_r}
    \partial \bar{w}^{j_1} \ldots \bar{w}^{j_r}}(z,w).
\EEQN
We clearly get
   $ M_1(F, G) =
    m \circ {\cal M}_1 (F \otimes G) = m \circ N (F \otimes G)$
and $M_r (F, G) = m \circ {\cal M}_r (F\otimes G)$. By induction
we find the following recursion formula
\BEQ {calMrRekur}
    {\cal M}_{r+1} = (N-r(n-r)-rH){\cal M}_r
\EEQ
where $H:=
z^i\frac{\partial}{\partial z^i}+
\bar{z}^i\frac{\partial}{\partial \bar{z}^i}+
w^i\frac{\partial}{\partial w^i}             +
\bar{w}^i\frac{\partial}{\partial \bar{w}^i}$. Note that this operator
vanishes on smooth complex-valued functions $F$ which are doubly homogeneous,
i.e. $F(\alpha z,\beta w)=F(z,w)$, and commutes with ${\cal M}_r$ and $N$.
This recursion can be solved and leads to the equation
\BEQN {calMreqn}
    {\cal M}_r = \prod_{s=0}^{r-1} (N - s(n-s)-sH)
\EEQN
Therefore all the operators ${\cal M}_r$ are
polynomials in the operator $N$ and $H$.

% ********** star-product **************

\section{The star-product on  $\protect\Complex \!\mbox{P}^n$
   and phase space reduction}

Before we construct the star-product on $\CP^n$ let us first
recall some well-known facts about the phase space
reduction of $\CNull{n+1}$ to $\CP^n$ with respect to the $U(1)$-action
$z\mapsto e^{{\bf i}\phi} z$, $\phi\in\Bbb{R}$. We basically use the
notation of \cite{abraham}.
An $ad^\ast$-equivariant momentum mapping is given by
\BEQ {DefJ}
    J(z) := -\frac{1}{2}x
\EEQ
and every $\mu \in \Bbb{R}^-$ is a regular value of $J$. Therefore
$J^{-1}(\mu)$ is a submanifold of $\CNull{n+1}$, the $(2n+1)$-sphere of
radius $\sqrt{-2\mu}$. The reduced phase space is then given by
$J^{-1}(\mu)/U(1) \cong \CP^n$. A symplectic form on $\CP^n$ is uniquely
defined by
\BEQ {CPnSympl}
    i_\mu^\ast \omega_0 = \pi_\mu^\ast \omega_\mu
\EEQ
where $i_\mu : J^{-1}(\mu) \to \CNull{n+1}$ is the inclusion map,
and
$\pi_\mu : J^{-1}(\mu) \to J^{-1}(\mu)/U(1)$ the projection on
the equivalence classes. For a $U(1)$-invariant function
$F\in \Cinfty$ we can define the reduced function
$F_\mu : \CP^n \to \Bbb{C}$ by
\BEQ {redfunc}
    F_\mu ([z]) := F \circ i_\mu (z)
    \qquad \mbox{ for } z \in J^{-1}(\mu)
\EEQ
Note that $F$ is $U(1)$-invariant iff $\{ F, J \} = 0$.
For two $U(1)$-invariant functions $F,G$ we have
\BEQ {poisson}
    \left( \{F, G\} \right)_\mu = \{ F_\mu, G_\mu \}_\mu
\EEQ
where $\{,\}_\mu$ denotes the Poisson bracket on $\CP^n$ which is
defined by $\omega_\mu$.

In the next Proposition we show how phase space reduction can be extended
to the star-products on $\CNull{n+1}$:

\begin{PROP} \label{algebraprop}
  Denote by
 ${\cal A}^0$ the subspace of $\cal A$ of those power series
  whose coefficients are $U(1)$-invariant.
  Suppose in addition that the power series $D$ \EQREF{theD} is equal to $1$.
  Then
  \begin{enumerate}
   \item A function $F$ in $\Cinfty$ is
        $U(1)$-invariant iff $F\ast J-J\ast F=0$ whence ${\cal A}^0$ is
        both a $\ast$- and a $\newstar$-subalgebra of $\cal A$.
   \item Fix a negative real number $\mu$. Denote by ${\cal J}_\mu$ and
        $\tilde{\cal J}_\mu$ the $\ast$-ideal and the $\newstar$-ideal
        of ${\cal A}^0$ generated by $J-\mu$, respectively.
        Then
        $\tilde{\cal J}_\mu$
        consists of those power series in ${\cal A}^0$ whose coefficients
        vanish on the $2n+1$-dimensional sphere $J^{-1}(\mu)$ and
        ${\cal J}_\mu$ is equal to $S^{-1}\tilde{\cal J}_\mu$.
   \item Denote by $\cal B$ (and by $\pi^*{\cal B}$) the space
        of formal power series with coefficients in
         $C^{\infty}(\CP^n)$ ($\pi^*(C^{\infty}(\CP^n)))$).
         Then we have the following direct sums:
        \BEQN{dirsum}
           {\cal A}^0=\pi^*{\cal B}\oplus {\cal J}_\mu
                     =\pi^*{\cal B}\oplus \tilde{\cal J}_\mu
        \EEQN
  \end{enumerate}
\end{PROP}
\begin{PROOF}
1. Using the definition of the Wick product \EQREF{WickDef} we see that
 $F\ast J-J\ast F = \frac{{\bf i}\lambda}{2}\{F,J\}$, and the latter
 vanishes iff $F$ is
 $U(1)$-invariant. In any associative algebra the subspace of all elements
commuting with a fixed one is a subalgebra. Since $S$ leaves ${\cal A}^0$
invariant (because $S$ commutes with the $U(1)$-action) it follows that
${\cal A}^0$ is also a $\newstar$-subalgebra of $\cal A$.

2. Let $g$ be a $U(1)$-invariant function vanishing on the
  $2n+1$-dimensional sphere $J^{-1}(\mu)$. Using the diffeomorphism
  $\CNull{n+1}\cong \Bbb{R}\times J^{-1}(\mu)$ and the $U(1)$-invariance
  of $g$ this function can be written as $g(z)=h(x,[z])$ where
  $h$ is a smooth complex-valued function on $\Bbb{R}\times \CP^n$.
  Then by Hadamard's trick
  \begin{eqnarray*} h(x,[z]) &=& h(x,[z])-h(-2\mu,[z])
             =\int_0^1\!dt\frac{d}{dt}h(tx+(1-t)(-2\mu),[z]) \\
  & = & -2
   (\int_0^1\!dt(\partial_1 h)(tx+(1-t)(-2\mu),[z]))(J(x)-\mu)
  \end{eqnarray*}
  which---by the obvious $U(1)$-invariance of the integral term---shows
  that $g$ is a $U(1)$-invariant multiple of $J-\mu$. But since
  $(J-\mu)g=(J-\mu)\newstar g=g\newstar(J-\mu)$ by \EQREF{Rnewstaru1}
  we obtain the asserted characterization of the $\newstar$-ideal
  $\tilde{\cal J}_\mu$. Moreover, because of the choice $D(x,\lambda)=1$
  and \EQREF {SxDx} we get
  \[ {\cal J}_\mu=(J-\mu)\ast{\cal A}^0=S^{-1}(S(J-\mu)\newstar S{\cal A}^0)
      =S^{-1}((J-\mu){\cal A}^0)= S^{-1}\tilde{\cal J}_\mu. \]

  3. Let $g$ be a smooth complex-valued $U(1)$-invariant function on
  $\CNull{n+1}$. Obviously $g=\pi^*g_\mu+(g-\pi^*g_\mu)$ where
  $g-\pi^*g_\mu$ clearly vanishes on $J^{-1}(\mu)$. Furthermore
  $\pi^*g$ vanishes on $J^{-1}(\mu)$ iff $g$ vanishes on $J^{-1}(\mu)$.
  This shows the first direct sum decomposition in \EQREF{dirsum}.
  Applying the linear bijection $S^{-1}$ to this first direct sum,
  observing that $Sf=f=S^{-1}f$ for all $f\in\pi^*{\cal B}$ and using part 2
  of this proposition we have shown the second direct sum in \EQREF{dirsum}.
\end{PROOF}

Mapping homogeneous functions to homogeneous functions, the bidifferential
operators $M_r$ (\EQREF{dieMs}) induce
corresponding bidifferential operators on $C^{\infty}(\CP^n)\times
C^{\infty}(\CP^n)$ by (recall that $x:=\bar{z}^iz^i$):
\BEQN {dietildeMs}
    \tilde{M}_r (\phi, \psi) ([z]) :=
             M_r(\phi\circ\pi,\psi\circ\pi)(z)
  =x^r\frac{\partial^r (\phi\circ\pi)}{\partial z^{i_1}\cdots
                    \partial z^{i_r}}(z)
                 \frac{\partial^r (\psi\circ\pi)}{\partial \bar{z}^{i_1}
                      \cdots \partial \bar{z}^{i_r}}(z).
\EEQN
where $\phi,\psi$ are two smooth complex-valued functions on $CP^n$.
In a completely analogous way, we may construct bidifferential operators
$\tilde{K}_r$ on $C^{\infty}(\CP^n)\times C^{\infty}(\CP^n)$ and
differential operators $\tilde{N}$
and $\tilde{\cal M}_r$ on $C^{\infty}(\CP^n\times\CP^n)$ from the
bidifferential
operators $K_r$ (\EQREF{homohomo}) on $\Cinfty\times \Cinfty$ and the
differential operators $N,{\cal M}_r$ (\EQREF{NOp}, \EQREF{calMR}) on
$C^{\infty}(\CNull{n+1}\times\CNull{n+1})$, respectively. The polynomial
relation \EQREF{calMreqn} also holds for $\tilde{N}$ and $\tilde{\cal M}_r$
after replacing $H$ by zero.

We can now state the main theorem:

\begin {THEOREM} \label {cpnstar}
    Let $\phi, \psi :\CP^n \to \Bbb{C}$ be two $C^\infty$-functions,
    let the bidifferential operators $\tilde{M}_r$ defined as is
    \EQREF{dietildeMs}, let $F$ and $G$ two smooth complex-valued
    $U(1)$-invariant functions on $\CNull{n+1}$, and $\mu$ a negative
    real number. Suppose that the series $D(x,\lambda)$ is equal to $1$.
   Then
    \begin{enumerate}
    \item
    \begin{eqnarray}
       \hspace{-3mm} \phi \mustar \psi \!\!
    &:= &\!\!\phi\psi+\sum_{r=1}^\infty \left(\frac{\lambda}{-2\mu}\right)^r
                   \sum_{s=1}^r \sum_{k=1}^s
                   \frac{1}{s!(s-k)!(k-1)!}
                     k^{r-1} (-1)^{r-k}~\tilde{M}_s(\phi ,\psi) \nonumber\\
                          & =:& \sum_{r=0}^\infty \lambda^r
                             \frac{1}{(-2\mu)^r} \tilde{K}_r (\phi, \psi).
                                       \label{KrCPn}
    \end{eqnarray}
    defines a star-product with first order commutator
       $ \frac{{\bf i}}{2} \{\phi, \psi\}_\mu$.
    \item
       Let $F, G$ be two $U(1)$-invariant functions
    then we have
    \BEQ {FGred}
        \left(F \newstar G \right)_\mu = F_\mu \mustar G_\mu.
    \EEQ
    \item
        There are the following isomorphisms of
        associative algebras:
    \BEQ {quotient}
          ({\cal B},\mustar) \cong
                 ({\cal A}^0,\newstar)/\tilde{\cal J}_\mu
                \cong ({\cal A}^0,\ast)/{\cal J}_\mu
    \EEQ
    \end{enumerate}
\end {THEOREM}
\begin{PROOF}
1. This can either be computed directly by using the fact that the
$\newstar$ product of a radial
function with a $U(1)$-invariant one is simply pointwise multiplication.
A more direct way of seeing this comes from the above proposition:
Since ${\cal B}\cong {\cal A}^0/\tilde{\cal J}_\mu$ as vector spaces
and $\tilde{\cal J}_\mu$ is a $\newstar$-ideal of ${\cal A}^0$ the
space $\cal B$ is equipped with an associative multiplication induced
by $\newstar$ where ``mod $\tilde{\cal J}_\mu$'' is translated into
``set $x$ equal to $-2\mu$''. Moreover, since $\ast$ and $\newstar$
are equivalent and the pointwise multiplication is commutative it is
easy to see that $\ast$ and $\newstar$ produce the same first order
commutator.

2. and 3. This is immediate from the above proposition and the first
part of this theorem since the linear map
$F\mapsto\pi^*F_\mu$ is nothing but the projection onto $\pi^*\cal B$
along the $\newstar$-ideal $\tilde{\cal J}_\mu$.
\end{PROOF}

\noindent
{\it Remarks:}
\begin{enumerate}
\item
In principle we could have obtained the star-product on $\CP^n$
directly from the Wick product by using the above second direct sum
decomposition of ${\cal A}^0$ into
$\pi^*{\cal B}$ and the $\ast$-ideal  ${\cal J}_\mu$: for two smooth
complex-valued functions $\phi$ and $\psi$ on $\CP^n$ we can form the Wick
product $\pi^*\phi\ast\pi^*\psi$ and recursively seperate off
order by order $\ast$-factors of $J-\mu$.
\item  For $D\not\equiv 1$
we still can construct a star-product on $\CP^n$ from the star-product
$\tilde{*}$ on $\CNull{n+1}$. Equation (\ref{thenewstar}) implies
that the star-product constructed in this way is obtained from
the star-product for $D \equiv 1$ by a ``reparametrisation''
of the formal parameter:
$\lambda \mapsto \lambda/(-2 \mu D(-2 \mu, \lambda))$. Such a
substitution is algebraically well defined because $\mu$ is
just a real number, and the formal power series $D(-2 \mu,\lambda)$
starts with 1, and obviously yields a star-product again.

\item An immediate consequence of this theorem  is the fact
that for $D\equiv 1$ the $\newstar$-commutator with $J$ coincides
with the classical Poisson-bracket and hence generates the same
$U(1)$-action. Thus, we may interpret $J$ for this choice of $D$ as a
``quantum momentum mapping''. In the general case, we can proceed
as follows: Using the facts  that the Wick-commutator of
  an arbitrary function $F$ and $J$ is
just  the Poisson bracket times $\frac{{\bf i}\lambda}{2}$ and has
no higher orders in $\lambda$, and that
the action of $\{ \cdot, J \}$ on $F$ is a derivation which commutes
with $\partial_x$ and therefore with $S$ and $S^{-1}$,
we conclude:
\BEQ {qmmap}
    F \newstar SJ - SJ \newstar F = \frac{{\bf i}\lambda}{2} \{F, J \}
\EEQ
This means that $SJ = D J $ is the quantum momentum mapping.
Hence, the case $D\equiv 1$ is distinguished by the fact that
the classical and quantum momentum maps coincide.

\end{enumerate}

\section{The domain $SU(1,n)/S(U(1)\times U(n))$}

In this section we shall sketch a modification
of the above results
to the noncompact dual (in the sense of Riemannian symmetric spaces
\cite{Hel 78}) of complex projective space: consider again $\CNull{n+1}$
as a complex manifold. Let $g^{ij}$ be a diagonal matrix with $(g^{ij})=
diag(-1,1,\ldots,1)$. Consider the
modified symplectic form
$\omega_1:=\frac{{\bf i}}{2}g^{ij}dz^i\wedge d\bar{z}^j$
and its corresponding Poisson bracket. The modified Wick product
$\ast_1$ is then given by
\BEQ {modWickDef}
    F \ast_1 G := \sum_{r=0}^\infty \frac{\lambda^r}{r!}
                g^{i_1j_1}\cdots g^{i_rj_r}
                \frac{\partial^r F}
                {\partial z^{i_1} \ldots \partial z^{i_r}}
                \frac{\partial^r G}
                {\partial \bar{z}^{j_1} \ldots \partial \bar{z}^{j_r}}.
\EEQ
This is easily seen to be an associative star-product corresponding to
the modified symplectic form $\omega_1$. Let $y(z):=g^{ij}z^i\bar{z}^j$.
Consider the open subset
$C^{n+1}_+\subset\CNull{n+1}$ defined by
\BEQ {dasC}
    C^{n+1}_+ := \{z\in \CNull{n+1}|y(z)>0\}.
\EEQ
It is clear that for each $z\in C^{n+1}_+$ the complex ray through
$z$ is also contained in $C^{n+1}_+$. Define
\BEQ{dasD}
   D^n := \pi(C^{n+1}_+)~\subset~\CP^n.
\EEQ
$D^n$ becomes a K\"ahler manifold with complex structure induced
by that of $\CP^n$ and symplectic form derived from the K\"ahler
potential $\log y$.
The noncompact semisimple Lie group $SU(1,n)$ acting canonically
on $\CNull{n+1}$ obviously preserves
the function $y$ and hence the open subset $C^{n+1}_+$. It is easy to
check that $SU(1,n)$ acts transitively on $D^n$ via its
action on the complex rays in $C^{n+1}_+$ thereby leaving the K\"ahler
structure of $D^n$ invariant. The isotropy subgroup of the standard
ray through $(1,0,\ldots,0)^T$ is clearly given by $S(U(1)\times U(n))$
hence $D^n$ is isomorphic to the irreducible hermitean symmetric space
$SU(1,n)/S(U(1)\times U(n))$ (see \cite{Hel 78}, p. 518).

It is not difficult to see that the whole construction mentioned in
the previous section can be transferred to $D^n$ by substituting
$\CNull{n+1}$
by $C^{n+1}_+$, homogeneous functions on $\CNull{n+1}$ by homogeneous
functions merely defined on $C^{n+1}_+$, the function $x$ by $y$,
radial functions by functions of $y$ and the bidifferential operators
$M_r$ by
\BEQ{diecheckMr}
    \check{M}_r(G,H):=
            y^rg^{i_1j_1}\cdots g^{i_rj_r}
    \frac{\partial^r G}{\partial z^{i_1}\cdots \partial z^{i_r}}
   \frac{\partial^r H}{\partial \bar{z}^{j_1}\cdots \partial \bar{z}^{j_r}}.
\EEQ

Observing that $J_1(z):=-\frac{1}{2}y(z)$ is a momentum map for the
canonical $U(1)$-action on $C^{n+1}_+$ corresponding to the symplectic
form $\omega_1$ we see that $D^n$ is a reduced phase space of $C^{n+1}_+$.
Consequently, both \PROPREF{algebraprop} and \THEOREMREF{cpnstar} (with
its explicit formula \EQREF{KrCPn}) are also
valid for the domain $D^n$ when modified along the lines given above.

%
% Vergleich mit Moreno
%

% ********** Vergleich *****************

\section{Comparison with existing results}

In this section we want to prove that in the case of $\CP^n \cong S^2$
our star-product $\mustar$ for $\mu=\frac{1}{2}$ is the same as the star
product constructed by Moreno and Ortega-Navarro in \cite {MO 83} iff
$D \equiv 1$.

To compare the two star-products we have to do some calculation in
the standard inhomogeneous complex co-ordinates
$v^i:=z^i/z^0, 1\leq i\leq n$ on $\CP^n$
First we see that the bidifferential operator $\tilde N$
can be written as ($(u^1,\ldots,u^n)$ being another set of inhomogeneous
co-ordinates on $\CP^n$):
\BEQ {NisLap}
    \tilde N = \left( 1+\sum_{j=1}^n v^j \bar{u}^j \right)
               \sum_{k,l=1}^n \left(v^k \bar{u}^l + \delta^{kl} \right)
               \frac{\partial}{\partial v^k}
               \frac{\partial}{\partial \bar{u}^l}.
\EEQ
For two locally defined real-analytic complex-valued functions
$\phi, \psi$ on $\CP^n$ we use their local holomorphic and
antiholomorphic continuation to obtain:
$(m \circ \tilde N (\phi, \psi)) (v, \bar{v}) =
     m \circ \Delta_{u\bar{u}} \phi (u, \bar{v}) \psi (v, \bar{u})$.
For these kind of functions we can inductively prove the following
local expression for the bidifferential operators $\tilde M_r$:
\BEQ {localMr}
    \tilde M_r (\phi, \psi) (v, \bar{v}) = m \circ
    \prod_{k=0}^{r-1} (\Delta_{u\bar{u}} + k(k-n))
    \phi (u, \bar{v}) \psi (v, \bar{u}) =:
           m\circ p_r(\Delta)_{u\bar{u}}\phi (u, \bar{v}) \psi (v, \bar{u})
\EEQ
Analogously we can conclude that the operators $\tilde K_r$
(\ref{KrCPn}) can again be thought of as polynomials in the Laplacian
which we denote by $\tilde{k}_r(\Delta)$.
We thus have to prove the recursion
formula \cite [Prop. 1.i] {MO 83} for these polynomials:
\BEQ {recKr}
    (r+1) \tilde{k}_{r+1}(\Delta) - \left[ \Delta \tilde k_r (\Delta)
    - \sum_{s=0}^{r-1} \frac{r! (r+3+s)}{(s+2)!(r-1-s)!} \tilde k_{r-s}
    (\Delta) \right]
    \stackrel {!}{=} 0
\EEQ

To prove this recursion formula we first write $\tilde{k}_r(\Delta)
=\sum_{t=1}^r\frac{1}{t!}p_t(\Delta)A^{(t)}_{r-t}$ (see \EQREF{Arsdef})
and get on the left hand side a linear combination of the
$p_t(\Delta)$. Noting that the terms of $p_1(\Delta)$ and $p_{r+1}(\Delta)$
automatically cancel we obtain for the coefficient of the $p_s(\Delta)$
for $r\ge 2$ and $s=2,\ldots, r$ the following expression:
\BEQN {coeff}
        A^{(s)}_{r+1-s} + \frac{1}{r+1} \sum_{t=s}^r
        \left[ {r+1 \choose t-1} + { r \choose t-1} \right]
        A^{(s)}_{t-s} - \frac{s}{r+1} \left(
        A^{(s-1)}_{r+1-s} - (s-1) A^{(s)}_{r-s} \right).
\EEQN
The polynomials $p_t(\Delta)$
are all linear independent therefore each
of the coefficients in (\ref{coeff}) has to vanish. In order to prove
that the expression (\ref{coeff}) vanishes
for $r\ge 2$ and $s=2, \ldots, r$ we
first note that in the case $s=r$ it can be checked directly. Using the
identity $A^{(s)}_{t+1-s} = A^{(s-1)}_{t+1-s} - s A^{(s)}_{t-s}$
for $s\ge1$, $t\ge s$ we can prove that (\ref{coeff}) vanishes by
induction on $r$.

\vspace{0.5cm}

\noindent
{\large\bf Acknowledgment}

\vspace{2mm}

\noindent The authors would like to thank E. Meinrenken for valuable
discussions.

\newpage

% ********** die Literatur *************

\end {document}